\documentclass[12pt]{article}
\usepackage{color}
\definecolor{darkblue}{rgb}{0.00,0.25,0.50}
\usepackage[colorlinks,filecolor=blue,citecolor=darkblue]{hyperref}

\usepackage{amsfonts,amssymb,amsmath,amsthm}
\usepackage{url}
\usepackage{enumerate}
\usepackage[ukrainian, russian, english]{babel}
\usepackage[cp1251]{inputenc}
\begin{document}\selectlanguage{ukrainian}
\thispagestyle{empty}

\title{}

\begin{center}
\textbf{\Large Оцінки найкращих $m$--членних тригонометричних наближень класів аналітичних функцій}
\end{center}
\vskip0.5cm
\begin{center}
А.~С.~Сердюк${}^1$, Т.~А.~Степанюк${}^2$\\ \emph{\small
${}^1$Інститут математики НАН
України, Київ\\
${}^2$Східноєвропейський національний університет імені Лесі
Українки, Луцьк\\}
\end{center}
\vskip0.5cm

%\address{Institute of Mathematics of the National Academy of
%Sciences of Ukraine\\ 3\\ Tereshenkivska st.\\ 01601\\ Kiev, Ukraine}

\begin{abstract}
 В метриках просторів $L_{s}, \ 1\leq s\leq\infty$, одержано  точні за порядком
оцінки знизу найкращих  $m$--членних тригонометричних наближень класів згорток періодичних функцій,  що
 належать одиничній кулі
простору $L_{p}, \ 1\leq p\leq\infty$, з твірним ядром ${\Psi_{\beta}(t)=\sum\limits_{k=1}^{\infty}\psi(k)\cos(kt-\frac{\beta\pi}{2})}$, $\beta\in \mathbb{R}$, коефіцієнти $\psi(k)$ якого  прямують до нуля не повільніше за геометричну прогресію.
Знайдені оцінки збіглися за порядком із наближеннями частинними сумами Фур'є  вказаних класів функцій в $L_{s}$--метриці, що дозволило також записати точні порядкові оцінки найкращих ортогональних тригонометричних наближень та тригонометричних поперечників зазначених класів.

\vskip 0.5cm

 In metric of spaces $L_{s}, \ 1\leq s\leq\infty$,
  we obtain exact in order estimates of best  $m$--term trigonometric approximations of classes of convolutions of periodic functions, that belong to unit ball of space  $L_{p}, \ 1\leq p\leq\infty$, with generated kernel $\Psi_{\beta}(t)=\sum\limits_{k=1}^{\infty}\psi(k)\cos(kt-\frac{\beta\pi}{2})$, $\beta\in \mathbb{R}$, whose coefficients $\psi(k)$ tend to zero not slower than  geometric progression. Obtained estimates coincide in order with approximation by Fourier sums of the given classes of functions in $L_{s}$--metric. This fact allows to write down exact order estimates of best orthogonal trigonometric approximation and trigonometric widths of given classes.

\end{abstract}

%%%%%%%%%%%%%%%%%%%%%%%%%%%%%%%%%%%%%%%%%%%%%%%%%%%%%%%%%%%%%%%%%%%%%%%%%

Нехай
$L_{p}$,
$1\leq p<\infty$, --- простір $2\pi$--періодичних сумовних в $p$--му
степені на $[0,2\pi)$  функцій $f:\mathbb{R}\rightarrow\mathbb{C}$ з нормою
$${\|f\|_{p}:=\Big(\int\limits_{0}^{2\pi}|f(t)|^{p}dt\Big)^{\frac{1}{p}}},$$
$L_{\infty}$ --- простір
$2\pi$--періодичних вимірних і суттєво обмежених  функцій $f:\mathbb{R}\rightarrow\mathbb{C}$ з нормою
$$\|f\|_{\infty}:=\mathop{\rm{ess}\sup}\limits_{t}|f(t)|.$$

Нехай $f:\mathbb{R}\rightarrow\mathbb{R}$ --- функція із $L_{1}$, ряд Фур'є  якої має вигляд
$$
\sum_{k=-\infty}^{\infty}\hat{f}(k)e^{ikx},
$$
де $\hat{f}(k):=\frac{1}{2\pi}\int\limits_{-\pi}^{\pi}f(t)e^{-ikt}dt,$
 $\psi(k)$ --- довільна фіксована послідовність дійсних чисел і $\beta$
--- фіксоване  дійсне число. Тоді якщо ряд
$$
\sum_{k\in \mathbb{Z}\setminus\{0\}}\frac{\hat{f}(k)}{\psi(|k|)}e^{i(kx+\frac{\beta\pi}{2} \mathrm{sign} k)}
$$
\noindent є рядом Фур'є деякої  функції $\varphi$ з $L_{1}$, то цю
функцію  називають $(\psi,\beta)$-похідною функції $f$ і позначають через
$f_{\beta}^{\psi}$ (див., наприклад, \cite[с. 132]{Stepanets1}).
Множину функцій $f$ у яких існує $(\psi,\beta)$-похідна
позначають через $L_{\beta}^{\psi}$.

Розглянемо одиничну кулю $B_{p}$ в просторі дійснозначних функцій з $L_{p}$, тобто множину функцій $\varphi:\mathbb{R}\rightarrow\mathbb{R}$ таких, що
$\|\varphi\|_{p}\leq1, \ 1\leq p\leq\infty$.
Якщо $f\in L^{\psi}_{\beta}$ і крім того $f^{\psi}_{\beta}\in
B_{p}$, то кажуть, що функція $f$ належить
класу
  $L^{\psi}_{\beta,p}$.

Підмножини неперервних функцій із $L^{\psi}_{\beta}$ та $L^{\psi}_{\beta,p}$ будемо позначати через $C^{\psi}_{\beta}$ та $C^{\psi}_{\beta,p}$ відповідно.

У випадку коли $\psi(k)=k^{-r},\ r>0$,  класи $L_{\beta,p}^{\psi}$, є
відомими класами
 Вейля-Надя $W_{\beta,p}^r$.

Послідовності $\psi (k),\ k\in \mathbb{N},$ що визначають класи
$L^{\psi}_{\beta,p}$, зручно розглядати як звуження на множину
натуральних чисел $\mathbb{N}$ деяких додатних, неперервних, опуклих
донизу функцій $\psi(t)$, $t\geq1$ таких, що $
\lim\limits_{t\rightarrow\infty}\psi(t)=0. $ Множину всіх таких
функцій $\psi(t)$ будемо позначати через ${\mathfrak M}$.

Наслідуючи О.І. Степанця (див., наприклад, \cite[с. 160]{Stepanets1}), кожній функції $\psi\in{\mathfrak M}$ поставимо у відповідність
характеристики
$$
\eta(t)=\eta(\psi;t)=\psi^{-1}\left(\psi(t)/2\right), \ \ \
\mu(t)=\mu(\psi;t)=\frac{t}{\eta(t)-t},
$$
де $\psi^{-1}$ --- обернена до $\psi$ функція i покладемо
$$
\mathfrak{M}^{+}_{\infty}=\left\{\psi\in \mathfrak{M}: \ \
\mu(\psi;t)\uparrow\infty, \ t\rightarrow\infty \right\}.
$$
Через ${\mathfrak M^{'}_{\infty}}$ позначимо  підмножину функцій
${\psi\in {\mathfrak M^{+}_{\infty}}}$ для яких величина
$\eta(\psi;t)-t$ обмежена зверху, тобто існує стала $K>0$ така,
що ${\eta(\psi;t)-t\leq K}, \ t\geq1$.

Як випливає з \cite[с. 1698]{Stepanets_Serdyuk_Shydlich} функції з множини $C^{\psi}_{\beta}$, де $\psi\in{\mathfrak M^{'}_{\infty}}$,
складаються із тих і тільки тих $2\pi$--періодичних функцій $f:\mathbb{R}\rightarrow\mathbb{R}$, які  допускають аналітичне продовження в смугу $|\mathrm{Im} z|\leq c, \ c>0$ комплексної площини. Отже, класи $C^{\psi}_{\beta,p}$ є класами
аналітичних функцій.

Природними представниками функцій з множини $\mathfrak{M}^{'}_{\infty}$  є
функції ${\psi(t)=\exp(-\alpha t^{r})}, \ {\alpha>0}, \ r\geq1$.

Нехай $f\in L_{s}$ і $\gamma_{m}$, $m\in\mathbb{N}$, --- довільний набір із $m$ цілих чисел. Величину
\begin{equation}\label{m_term1}
  e_{m}(f)_{s}=\inf\limits_{\gamma_{m}}\inf\limits_{c_{k}\in\mathbb{C}}\|f(x)-\sum\limits_{k\in \gamma_{m}}c_{k}e^{ikx} \|_{s}, \  \ 1\leq s\leq\infty,
\end{equation}
називають найкращим $m$--членним тригонометричним наближенням функції ${f}$ в метриці простору $L_{s}$.
В більш загальній ситуації величини виду (\ref{m_term1}) при $s=2$
були введені С.Б. Стєчкіним \cite{Stechkin}
з метою встановлення критерію абсолютної збіжності ортогональних рядів.

Для довільного класу $F$ із $L_{s}$ покладемо
\begin{equation}\label{m_term}
  e_{m}(F)_{s}:=\sup\limits_{f\in F}e_{m}(f)_{s},  \  \ 1\leq s\leq\infty.
\end{equation}

 Порядки спадання до нуля при $n\rightarrow\infty$ величин (\ref{m_term}) при   $F=L_{\beta,p}^{\psi}$ для різних співвідношень між параметрами $p$ i $s$ за умови $\psi\in B$, де $B$ --- множина
незростаючих додатних функцій $\psi(t)$, $t\geq 1$, для кожної з
яких існує додатня стала $K$ така, що $
\frac{\psi(t)}{\psi(2t)}\leq K, \ \  t\geq 1 $,  та при деяких додаткових  умовах  на функцію $\psi$ досліджувались у роботах \cite{Fedorenko2000}--\cite{Fedorenko2004}.

В даній роботі розглядається задача про знаходження точних порядкових оцінок величин $e_{m}(L_{\beta,p}^{\psi})_{s}$, $1\leq p,s\leq\infty$, $\beta\in\mathbb{R}$ у випадку коли $\psi\in{\mathfrak M^{'}_{\infty}}$.

Окрім величин $ e_{m}(F)_{s}$ в роботі для класів $F=L_{\beta,p}^{\psi}$ розглядаються   величини вигляду
$$
 e^{\bot}_{m}(F)_{s}=\sup\limits_{f\in F}\inf\limits_{\gamma_{m}}\|f(x)-\sum\limits_{k\in \gamma_{m}}\hat{f}(k)e^{ikx} \|_{s}, \ \ 1\leq s\leq\infty, \ m\in\mathbb{N},
$$
які називають найкращими ортогональними тригонометричними наближеннями класу $F=L_{\beta,p}^{\psi}\subset L_{s}$ в метриці простору $L_{s}$,
  а також величини
 $$
  d^{\top}_{m}(F)_{s}:=\inf\limits_{\gamma_{m}}\sup\limits_{f\in F}\inf\limits_{c_{k}\in\mathbb{C}}\|f(x)-\sum\limits_{k\in \gamma_{m}}c_{k}e^{ikx} \|_{s}, \  \  \ 1\leq s\leq\infty,
$$
які називають
 тригонометричними поперечниками класу ${F}$ в метриці простору $L_{s}$.

Позначимо через ${\cal E}_{n}(F)_{s}$ наближення сумами Фур'є класу $F\subset L_{s}$ в метриці простору $L_{s}$,
тобто величини
  вигляду
  $$
  {\cal E}_{n}(F)_{s}=\sup\limits_{f\in
F}\|f(x)-\sum\limits_{k=-n+1}^{n-1}\hat{f}(k)e^{ikx}\|_{s}, \ 1\leq s\leq\infty.
  $$

З означень величин $ e_{m}(F)_{s}$, $ e^{\bot}_{m}(F)_{s}$, $ d^{\top}_{m}(F)_{s}$ i ${\cal E}_{n}(F)_{s}$
випливає, що при ${m=2n, 2n-1}, \ n\in\mathbb{N}$ мають місце нерівності
\begin{equation}\label{ineq}
 e_{m}(F)_{s}\leq
                                     \begin{array}{c}
                                      e_{m}^{\bot}(F)_{s} \\
                                        d^{\top}_{m}(F)_{s} \\
                                     \end{array}
 \leq {\cal E}_{n}(F)_{s}.
\end{equation}

Зазначимо, що величини $ e_{m}(F)_{s}$, $ e^{\bot}_{m}(F)_{s}$, $ d^{\top}_{m}(F)_{s}$ і ${\cal E}_{n}(F)_{s}$ для різноманітних класів функцій $F$ як однієї, так і багатьох змінних вивчались багатьма авторами. З детальною історією досліджень цих величин та відповідною бібліографією можна ознайомитись, наприклад, в роботах \cite{Romanyuk2003}--\cite{Serdyuk_grabova}.

{\bf Теорема 1.} {\it Нехай  $\psi\in{\mathfrak M^{'}_{\infty}}$, $1\leq p,s\leq\infty$, $\beta\in\mathbb{R}$.
 Тоді  мають місце порядкові оцінки}
\begin{equation}\label{theorem}
  e_{2n}(L_{\beta,p}^{\psi})_{s}\asymp  e_{2n-1}(L_{\beta,p}^{\psi})_{s}\asymp \psi(n).
\end{equation}

{\bf Доведення теореми.}
В силу теореми 6.8.2 роботи \cite[с. 48]{Stepanets2}, якщо $\psi\in{\mathfrak M^{'}_{\infty}}$, то
\begin{equation}\label{Stepanets}
 {\cal E}_{n}(L^{\psi}_{\beta,p})_{s} \asymp \psi(n), \ \ 1\leq p,s\leq\infty, \ \ \beta\in\mathbb{R}.
\end{equation}
Згідно зі співвідношеннями (\ref{ineq}) i (\ref{Stepanets}) одержуємо
$$
  e_{2n}(L_{\beta,p}^{\psi})_{s}\leq e_{2n-1}(L_{\beta,p}^{\psi})_{s}\leq C_{0}\psi(n), \ \ 1\leq p,s\leq\infty,
$$
де $C_{0}$ --- деяка додатня стала.
Знайдемо відповідну оцінку знизу для величини $e_{2n}(L_{\beta,p}^{\psi})_{s}$. Доозначимо послідовність $\psi(k)$  у точці $k=0$ за допомогою рівності $\psi(0)=\psi(1)$.
 Розглянемо функцію
$$
f^{*}(t)=f^{*}(\psi;n;t):=C_{1}\Big(\frac{\psi(1)}{2(n+A)^{2}}+ \sum\limits_{k=1}^{n}\frac{\psi(k)}{(n-k+A)^{2}}\cos kt\Big),
$$
де $C_{1}$ та $A$ --- деякі додатні сталі, які будуть визначені пізніше.

Оскільки
$$
\|(f^{*})_{\beta}^{\psi}(t)\|_{p}=C_{1}\Big\|\sum\limits_{k=1}^{n}\frac{1}{(n-k+A)^{2}}\cos \Big(kt+\frac{\beta\pi}{2}\Big)\Big\|_{p}\leq
$$
$$
\leq
2\pi C_{1}\sum\limits_{k=1}^{n}\frac{1}{(n-k+A)^{2}},
$$
то очевидно, що вибравши сталу $C_{1}$ так, щоб $2\pi C_{1}\sum\limits_{k=1}^{n}\frac{1}{(n-k+A)^{2}}\leq1$, отримаємо включення $f^{*}\in L_{\beta,p}^{\psi}$.
Покажемо, що
\begin{equation}\label{eq5}
e_{2n}(f^{*})_{s}\geq C_{2}\psi(n), \ n\in\mathbb{N},
\end{equation}
де $C_{2}$ --- деяка додатня стала.

З цією метою скористаємось співвідношенням двоїстості (див., наприклад, \cite[с. 42]{Korn})
\begin{equation}\label{eq1}
 e_{m}(f)_{s}=\inf\limits_{\gamma_{m}}{\mathop{\sup}\limits_{
 h\in L^{\bot}(\gamma_{m}),\atop \|h\|_{s'}\leq1
 }}\bigg|\int\limits_{-\pi}^{\pi}f(t)h(t)dt \bigg|, \ m\in\mathbb{N},
\end{equation}
де $\frac{1}{s}+\frac{1}{s'}=1$, a запис
$h\in L^{\bot}(\gamma_{m})$ означає, що
$$
\int\limits_{-\pi}^{\pi}h(t)e^{ikt}dt=0, \ \ k\in\gamma_{m}.
$$

Для довільного набору $\gamma_{2n}$ із $2n$  цілих чисел візьмемо довільне ціле число $k^{*}$ таке, що  $k^{*}\in[-n,n]$ і $k^{*}\notin \gamma_{2n}$.
Покладемо
$$
T(t):=\frac{1}{2\pi}e^{-i k^{*}t}.
$$

Очевидно, що  $T\in L^{\bot}(\gamma_{2n})$ i $\|T\|_{s'}\leq1, \ 1\leq s\leq\infty, \ \frac{1}{s}+\frac{1}{s'}=1$,
а отже в силу співвідношення  (\ref{eq1}) та рівності
 $$
  \int\limits_{-\pi}^{\pi}e^{ikt}e^{imt}dt=
{\left\{\begin{array}{cc}
0, \ & k+m\neq 0, \\
2\pi, &
k+m=0, \
  \end{array} \right.}  \ \ k,m\in\mathbb{Z},
$$
 отримуємо оцінку
$$
 e_{2n}(f^{*})_{s}=\inf\limits_{\gamma_{2n}}{\mathop{\sup}\limits_{
 h\in L^{\bot}(\gamma_{2n}),\atop \|h\|_{s'}\leq1
 }}\bigg|\int\limits_{-\pi}^{\pi}f^{*}(t)h(t)dt \bigg|\geq\inf\limits_{\gamma_{2n}}
\bigg|\int\limits_{-\pi}^{\pi}f^{*}(t)T(t)dt \bigg|=
$$
$$
=\frac{C_{1}}{4\pi}\inf\limits_{\gamma_{2n}}
\bigg|\int\limits_{-\pi}^{\pi}\sum\limits_{|k|\leq n}\frac{\psi(|k|)}{(n-|k|+A)^{2}}e^{-i kt}e^{-i k^{*}t}dt \bigg|
=\frac{C_{1}}{2}\inf\limits_{\gamma_{2n}}
\frac{\psi(|k^{*}|)}{(n-|k^{*}|+A)^{2}}\geq
$$
\begin{equation}\label{eq2}
 \geq\frac{C_{1}}{2}\min\limits_{0\leq k\leq n}\frac{\psi(k)}{(n-k+A)^{2}}=\frac{C_{1}}{2}\min\limits_{1\leq k\leq n}\frac{\psi(k)}{(n-k+A)^{2}}.
\end{equation}
Покажемо, що функція ${\psi_{n}(t)=\frac{\psi(t)}{(n-t+A)^{2}}}$ при певному виборі сталої $A$ не зростає на проміжку $[1,n]$.
Легко бачити, що
$$
\psi'_{n}(t)=\bigg(\frac{\psi(t)}{(n-t+A)^{2}}\bigg)'=\frac{2\psi(t)}{(n-t+A)^{3}}+\frac{\psi'(t)}{(n-t+A)^{2}}=
$$
\begin{equation}\label{eq3}
 = \frac{\psi(t)}{(n-t+A)^{3}}\Big(2-\frac{|\psi'(t)|}{\psi(t)}(n-t+A)\Big), \ \ \psi'(t):=\psi'(t+0).
\end{equation}

Далі скористаємось тим, що для $\psi\in\mathfrak{M}_{\infty}^{+}$ за умови $\mu(t)\geq b>0$ має місце нерівність (див. \cite[с. 1251]{Serdyuk_Stepaniuk2014})
\begin{equation}\label{lemma}
\frac{\psi(t)}{|\psi'(t)|}\leq4\big(1+\frac{1}{b}\big)\left(\eta(t)-t\right),
\ \ t\geq 1.
\end{equation}

Оскільки $\psi\in\mathfrak{M'}_{\infty}$, то існує стала $K_{0}>0$, така, що $\eta(t)-t\leq K_{0}$, а отже
\begin{equation}\label{f1}
  \mu(t)=\frac{t}{\eta(t)-t}\geq\frac{1}{\eta(t)-t}\geq \frac{1}{K_{0}}
\end{equation}
i застосовуючи (\ref{lemma}) при $b=\frac{1}{K_{0}}$, маємо
\begin{equation}\label{f2}
  \frac{|\psi'(t)|}{\psi(t)}\geq\frac{1}{4(K_{0}+1)(\eta(t)-t)}\geq\frac{1}{4(K_{0}+1)K_{0}}.
\end{equation}

Враховуючи (\ref{f2}), отримуємо
 \begin{equation}\label{eq4}
 2-\frac{|\psi'(t)|}{\psi(t)}(n-t+A)\leq 2-\frac{1}{4(K_{0}+1)K_{0}}(n-t+A)\leq2-\frac{A}{4(K_{0}+1)K_{0}}.
\end{equation}

З (\ref{eq3}) i (\ref{eq4}) випливає, що при ${A\geq8K_{0}(K_{0}+1)}$ справедлива нерівність ${\psi'_{n}(t)\leq0}, \ t\geq1$,
тобто функція  $\psi_{n}(t)$ не зростає.
Тому
\begin{equation}\label{f}
\min\limits_{1\leq k\leq n}\frac{\psi(k)}{(n-k+A)^{2}}=\frac{\psi(n)}{A^{2}}.
\end{equation}
З (\ref{eq2}) та (\ref{f}) отримуємо (\ref{eq5}).
 Теорему доведено.

 Теорема 1 разом зі  співвідношеннями (\ref{ineq}) та (\ref{Stepanets}) дозволяють записати наступне твердження.

{\bf Теорема 2.} {\it Нехай  $\psi\in{\mathfrak M^{'}_{\infty}}$, $1\leq p,s\leq\infty$, $\beta\in\mathbb{R}$  і $m\in\mathbb{N}$.
 Тоді
 $$
 e_{m}(L_{\beta,p}^{\psi})_{s}\asymp  e_{m}^{\bot}(L_{\beta,p}^{\psi})_{s}\asymp  d_{m}^{\top}(L_{\beta,p}^{\psi})_{s}\asymp \psi\Big(\Big[\frac{m+1}{2}\Big]\Big),
 $$
де запис $[a]$ означає цілу частину дійсного числа $a$.}

E-mail: \href{mailto:serdyuk@imath.kiev.ua}{serdyuk@imath.kiev.ua},
\href{mailto:tania_stepaniuk@ukr.net}{tania$_{-}$stepaniuk@ukr.net}

\end{document}